\documentclass[11pt, a4paper]{article}

\usepackage{amsmath,amssymb, amsthm}
\usepackage{mathtools}
\usepackage[utf8]{inputenc}
\usepackage{mathrsfs}

\usepackage{fullpage}

\usepackage{authblk}
    
\usepackage{dsfont}
\usepackage{bbm}
    
\usepackage{graphicx}
\graphicspath{ {images/} }

\usepackage{wrapfig}

\usepackage{verbatim}

\usepackage{hyperref}

\usepackage{xcolor}
\hypersetup{
  colorlinks   = true, 
  urlcolor     = {blue!90!black}, 
  linkcolor    = {blue!90!black}, 
  citecolor   = {red!90!black} 
}

\newcommand{\R}{\mathbb{R}}
\newcommand{\N}{\mathbb{N}}
\newcommand{\Z}{\mathbb{Z}}
\renewcommand{\d}{{\mathrm{d}}}

\renewcommand{\P}{\mathbb{P}}
\newcommand{\E}{\mathbb{E}}
\newcommand{\PP}[1]{\mathbb{P}\left\{#1\right\}}
\newcommand{\EE}[1]{\mathbb{E}\left\{#1\right\}}
\newcommand{\Ec}[2]{\mathbb{E}\left[#1 \middle| #2 \right]}

\renewcommand{\Pr}{{\mathcal{P}}}
\newcommand{\U}{{\mathcal{U}}}

\newcommand{\1}{\mathds{1}}

\definecolor{ViktorColor}{rgb}{0.47, 0.00, 0.97}

\theoremstyle{plain}
\newtheorem{thm}{Theorem}[section]

\newtheorem{cor}[thm]{Corollary}
\newtheorem{rmk}[thm]{Remark}

\theoremstyle{definition}

\author[1,2]{Viktor Bezborodov \thanks{Email: \texttt{viktor.bezborodov@pwr.edu.pl}}} 

\affil[1]{
	{Wroc\l{}aw University of Science and Technology,
		 Faculty of Electronics 
	 }}
\affil[2]{
	{University of Goettingen,  Institute for Mathematical Stochastics
}}

\title{Non-triviality in a totally asymmetric
one-dimensional Boolean
percolation
 model on a half-line}

\begin{document}

\maketitle

\begin{abstract}
 
It is well known that there are two regimes   in a standard  one-dimensional Boolean percolation model: either
 the entire space is covered a.s.,
 or 
 the covered volume fraction is strictly less than one.
 The aim of this work is to demonstrate
 that
 there is a third possibility
  in 
   a Boolean model
   with totally asymmetric grains
   on a half-line: 
   a.s.
  there is no unbounded component,
  but the covered volume fraction is one. An explicit condition
  	is given
  	characterizing  the existence of an unbounded occupied component.
 
\end{abstract}


\section{Introduction}

Let
$\R_+ = [0,\infty )$ and let 
 $\Pr$ be a  Poisson point process on $\R \times (0,+\infty)$
with intensity measure
$\lambda du \times \mu$.
Here $\mu$ is  a probability distribution on $(0,+\infty)$,
and $\lambda > 0$.
We consider  totally asymmetric 
Boolean percolation (TABP) on $\R$ and $\R_+$ defined as follows.
On $\R$
the occupied 
part of the space
is the union  
$$\U =  \bigcup\limits _{(u, \rho _u) \in \Pr} [u, u+\rho _u] .$$
On $\R_+$
the occupied part of the space is
 $${\U _+ =  \bigcup\limits _{(u, \rho _u) \in \Pr: u \geq 0} [u, u+\rho _u] }.$$
The model  just described
belongs to the class
of the germ-grain models \cite[Chapter 6]{MeckeBook13}, \cite[Chapter 4]{SW08}.
The intervals $[u, u +\rho _u]$,  $(u, \rho _u) \in \Pr$,
are the grains, and the points $u$
are the germs.
A discrete-space equivalent of TABP appears in \cite[Section 3]{Lamp70},
where the grains are interpreted as fountains
that wet sites to their right. 

The properties of 
the corresponding 
symmetric model are well known. 
Let $\ell$ be one-dimensional Lebesgue measure. 
For $Q = \R _+$ or $Q = \R$, 
define the covered volume fraction (\cite[Chapter 5]{MR96}) as
\begin{equation}
\lim\limits _{n \to \infty} \frac{\ell \left([-n, n] \cap Q \cap \U \right)}{\ell \left([-n, n] \cap Q  \right)}
\end{equation}
Let $\rho$ be a random variable
with distribution $\mu$.
In the standard one-dimensional Boolean model
 the grains are given by 
$[u - \rho _u, u +\rho _u]$, 
and only two regimes exist (\cite[Sections 3.2 and 5.1]{MR96}):

\begin{itemize}
	\item If $\E \rho < \infty$
	 a.s.
	no unbounded occupied component occurs and	
	 the covered volume fraction is $1 - e^{-2 \lambda \E \rho} \in (0,1)$.

	\item If $\E \rho = \infty$
		for any $\lambda > 0$
	 a.s. the entire space is covered.	
\end{itemize}
These two  regimes
exhaust all possible scenarios
also
for  TABP on $\R$.
Indeed, 
$$\U =  \bigcup\limits _{(u, \rho _u) \in \Pr} [u, u+\rho _u] = \bigcup\limits _{(u, \rho _u) \in \Pr} \Big[ u+\frac {\rho _u}{2} - \frac {\rho _u}{2}, u+\frac {\rho _u}{2} + \frac {\rho _u}{2}\Big] 
= \bigcup\limits _{(v, \varrho _v) \in \Pr'} [ v - \varrho _v, v + \varrho _v],
$$
where $\Pr' = \{ (v, \varrho _v): v = u + \frac {\rho _u}{2}, \varrho _v = \frac {\rho _u}{2} \text{ for some } (u, \rho _u) \in \Pr \} $
is a Poisson point process on $\R \times (0,+\infty)$ with intensity $\lambda du \times \mu _2$, where $\mu _2$
is the distribution of $\rho/2$. Note that $\Pr'$ 
is the image of $\Pr$ under the linear transformation $(a,b) \mapsto \left(a+\frac b2, \frac b2\right)$.
 Thus we are essentially in the settings
of the Poisson Boolean model \cite{MR96}. Hence
for TABP on $\R$ given by $\U$

\begin{itemize}
	\item If $\E \rho < \infty$
	a.s.
	no unbounded occupied component occurs
    and	
	the covered volume fraction is $1 - e^{- \lambda \E \rho} \in (0,1)$.

	\item If $\E \rho = \infty$
	for any $\lambda > 0$
	a.s. the entire space is covered.	
\end{itemize}

The aim of the this work is to demonstrate that
for TABP on $\R _+$ (given by $\U_+$)
there is another alternative.
It appears in item (iii)
of Theorem \ref{main thm}, where $\E \rho = \infty$ and
the covered volume fraction is $1$ and despite that no unbounded component exists.
Furthemore, the necessary and sufficient conditions for the existence of an unbounded occupied component
for $\U$ and $\U_+$ are different.

\begin{thm} \label{main thm}
	Consider  totally asymmetric Boolean percolation on $\R _+$.
	
	 \begin{itemize}
	 	\item [] (i) Assume that $\E \rho < \infty$.
	 	Then   a.s. no unbounded occupied component exists, 
	 	and the covered volume fraction 
	 	is $1- e ^{- \lambda \E \rho}$.

	 	\item [] (ii) Assume that $\int\limits _0 ^\infty
	 	e^{-  {{ \lambda}} \E( \rho \wedge t)} dt < \infty$. 
	 	Then a.s. there exists an unbounded occupied component. 
	 	
	 		 	\item [] (iii) Assume that
	 		 	$\E \rho = \infty$ 
	 		 	and
	 		 	 $\int\limits _0 ^\infty
	 	e^{- {{ \lambda}} \E( \rho \wedge t)} dt = \infty$. 
	 	Then a.s. there is no unbounded occupied component,
	 	however the expected size of the occupied component
	 	containing $1$
	 	is infinite
	 	and the covered volume fraction is $1$.

	 \end{itemize}
\end{thm}

\begin{cor}
	Consider  totally asymmetric Boolean percolation on $\R _+$.
	A.s.
	there exists an unbounded occupied component if
	and only if $\int\limits _0 ^\infty
	e^{-  { \lambda} \E( \rho \wedge t)} dt < \infty$.
	A.s. there is no unbounded occupied component 
	if and only if $\int\limits _0 ^\infty
	e^{-  { \lambda} \E( \rho \wedge t)} dt = \infty$.
\end{cor}

	Note that 
	the assumption $\int\limits _0 ^\infty
	e^{-  {{ \lambda}} \E( \rho \wedge t)} dt < \infty$ in (ii) of Theorem \ref{main thm}
   implies 
		$\E \rho = \infty$. All possible distributions $\mu$
		are covered by Theorem \ref{main thm} and items (i), (ii), and (iii) constitute  distinct 
		cases.
It is interesting to note that in the settings of (iii) 
the expected size of an occupied component is infinite, yet despite that
a.s.
every occupied component is finite.


In \cite{frogL} similar  results are collected for a discrete version of the model. 
In particular, the discrete space equivalent of Theorem \ref{main thm}
can be found in \cite[Section 3]{frogL}.
As one might expect,
	 in
	dimension $\d \geq 2$
the geometry    of the  Boolean model 
and the interplay between occupied and vacant components
  are more intricate. There is a 
  substantial amount of literature on the subject.
  We mention only a few recent works 
on the connectivity properties (such as 
the decay of the probability of reaching a remote point from the origin
in the subcritical regime, or the sharpness of the phase transition) \cite{DCRT20, ATT18}  
and the capacity functional (that is, the probability that the occupied part
of the space intersects a given compact set) \cite{CapBool}.

The proof of Theorem \ref{main thm}
is located Section \ref{sec 2}.
An example of a  distribution satisfying conditions of 
 Theorem \ref{main thm}, (iii),
 can be found in Remark \ref{trepidation}.

\section{Proof of Theorem \ref{main thm}}\label{sec 2}

By \cite[Theorem 9.3.5]{SW08}, 
the Boolean model is ergodic,
thus the sequence 
$\{\psi _n\}_{n \in \Z}$
defined by
\begin{equation}
\psi _n =  \ell \left([n,n+1] \cap \U \right) 
\end{equation}
is ergodic too. 
Hence 
\begin{equation}\label{tawdry}
\lim\limits _{x \to +\infty} \frac{ \ell \left([0,x] \cap\U \right) }{x}
= \E \psi _0 = \PP{ 0 \in \U }
= 1 -  e ^{- \lambda \E \rho}.
\end{equation}

\textbf{Proof of Theorem \ref{main thm}.}
We start with (i).
 Assume  $\E \rho < \infty$.
The fact that a.s. no unbounded occupied component exists
is the content of \cite[Theorem 3.1]{MR96}.
We have 
	\begin{multline*}
		 \E \# \big\{ (u, \rho _u) \in \Pr: u < 0,  (u + \rho _u) \geq 0  \big\}
		 = \E \# \big\{  \Pr \cap  \big\{  (u, p): u < 0, p \geq |u|   \big\}  \big\}
		 \\
		 =  \iint \limits _{ \substack{ (u, p): u < 0, \\ p \geq |u|    }} 
		 \lambda du \mu (dp)  
		 = \lambda \int \limits _{p \geq 0  } 
		  \mu (dp)  \int \limits _{-p \leq u \leq 0 } du
		  = \lambda \int  \limits _{p \geq 0  } p \mu (dp)  = \lambda \E \rho < \infty. 
	\end{multline*}
Hence a.s. $ \# \big\{ (u, \rho _u) \in \Pr: u < 0,  (u + \rho _u) \geq 0  \big\} < \infty$
and
\begin{equation}
 \sup\limits_{(u, \rho _u) \in \Pr: u < 0  } (u + \rho _u)  < \infty. 
\end{equation}
Consequently by  \eqref{tawdry} a.s.
\begin{equation}
\lim\limits _{x \to +\infty} \frac{ \ell \left([0,x] \cap \U_+ \right) }{x}
= \lim\limits _{x \to +\infty} \frac{ \ell \left([0,x] \cap\U \right) }{x}
=
1 -  e ^{- \lambda \E \rho}.
\end{equation}

Now we proceed to (ii) and (iii). 
For $t \geq 0$ the probability that $t$ is not covered 
by $\U _+$
\begin{align} \label{zaniness}
	\PP{ t \notin \U_+} & = \PP{  \Pr \cap \big\{  (u, p): 0 \leq u \leq t, p \geq t - u   \big\} = \varnothing  }
	\\
	& 
	=
	\exp\left\{ - \iint \limits _{(u, p): 0 \leq u \leq t, p \geq t - u  } 
	\lambda du \mu (dp)    \right\} \notag 
	\\
	& 
	=
	\exp\left\{ -  \lambda \int \limits _{0 } ^ \infty 
	\mu (dp) 
	\int \limits _{(t-p) \vee 0 } ^ t  du    \right\} \notag  
	=
	\exp\left\{ -  \lambda \int \limits _{0 } ^ \infty 
	(t \wedge p)
	\mu (dp) 
	\right\}
	\\
	&
	= e^{-  \lambda \E (t \wedge \rho) }.
\end{align}

Set $\tau_1=0$, $\xi_ 1 = \inf \{t > \tau_1: t \in \U _+  \} $, $\tau _2 = \inf\{ t > \xi_1: t \notin \U _+ \}$, 
 and so on, so that for $n \in \N$, $n \geq 2$
\[
\tau _n = \inf\{ t > \xi_{n-1}: t \notin \U _+ \}, \ \ \ \xi_ n = \inf \{t > \tau _n: t \in \U _+  \} 
\]
A.s. on $\{\tau _n < \infty\}$, $[\tau_{n}, \xi_n ]$ is the $n$-th connected component of $\R _+ \setminus \U_+$,
whereas $[ \xi_n, \tau_{n+1} ]$ is the $n$-th connected component of $ \U_+$.
A.s. on $\{\tau _n < \infty\}$, 
$$\xi _n = \min \{u > \tau _n:  (u, \rho _u) \in \Pr \text{ for some } \rho _u \in \R_+ \},$$
therefore given  $\{\tau _n < \infty\}$, $\xi _1 - \tau _1, \xi _2 - \tau _2, ...,\xi _n - \tau _n$ 
is a sequence of i.i.d random variables distributed 
exponentially with mean $ 1/\lambda$.

Let $N_v$ be the total number of vacant components, $N_v = \min\{ n \in \N : \tau _{n+1} = \infty  \}$ 
($N_v$
may take the value $+\infty$).
We have
\begin{equation}\label{goon}
	\E \ell (\R _+ \setminus \U_+) = \E \Ec{ \ell (\R _+ \setminus \U_+)}{N_v} =
   \E	\Ec{ \sum\limits _{i=1}^{N_v}(\xi _i - \tau _i)}{N_v}
	\\
	= \E \frac{N_v}{\lambda} =  \frac{\E N_v}{\lambda}.
\end{equation} 
On the other hand by \eqref{zaniness}
\begin{equation}\label{surly = sullen, arrogant}
	\E \ell (\R _+ \setminus \U_+) = \int _0 ^ \infty  \PP{ t \notin \U_+} dt=
	\int _0 ^ \infty  e^{-  \lambda \E (t \wedge \rho) } dt.
\end{equation}

Conditionally on $\{ \tau _n <\infty \}$, 
the distribution of
$\tau_{n+1}-\xi_n$ is the same
as the (unconditional) distribution of $\tau_2-\xi_1$.
Therefore
 $N_v$ has a geometric distribution with parameter ${p = \P \{ \tau _2 = \infty \}}$,
that is, ${\P \{ N_v = m \} = p (1-p)^{m-1}}$, $m = 1,2,...$;
if $p = 0$, then a.s. $N_v = \infty$.
In particular,  $\E N_v = 1/p$ if $p >  0$.
Note that $\E N_v < \infty$ if and only if $p> 0$.
By \eqref{goon} and \eqref{surly = sullen, arrogant}
\begin{equation} \label{grease}
	 \E N_v = \lambda \int _0 ^ \infty  e^{-  \lambda \E (t \wedge \rho) } dt.
\end{equation}

Assume that $\int\limits _0 ^\infty
	e^{-  {{ \lambda}} \E( \rho \wedge t)} dt < \infty$. Then 
	by \eqref{grease}, $\E N_v < \infty$
	and hence $N_v$ is a.s. finite. In this case the $N_v$-th
	occupied component is unbounded, and (ii) is proven.

	Assume that
	$\E \rho = \infty$ 
	and
	$\int\limits _0 ^\infty
	e^{- {{ \lambda}} \E( \rho \wedge t)} dt = \infty$. 
	By \eqref{grease}, $\E N_v = \infty$, hence $p = 0$
	and $N_v = \infty $ a.s. Thus, there are a.s. infinitely many 
	vacant components, that is there is no unbounded occupied component. 
	The expected size of an occupied  component 
	exceeds $\E\rho$, hence it is infinite.
	Let $l_1$, $l_2$, $...$ be the lengths
	of the consecutive occupied components (that is, $l_m = \tau _{m+1} - \xi _m$).
	Then $\E l _i = \infty$, $i = 1,2,...$,
	and
	a.s. $\frac 1n \sum\limits_{i=1}^n l_i \to \infty$, $n \to \infty$.
	Since a.s.
	the average length of the vacant component
	 $\frac 1n \sum\limits_{i=1}^n (\xi _i - \tau _i) \to \frac 1 \lambda$,
	the covered volume fraction is 1. 
	\qed

\begin{rmk}\label{trepidation}
	An example of $\lambda$ and $\rho $ satisfying 
	conditions of item (iii)
	of Theorem \ref{main thm}
	is $\lambda = 1$ and
	\begin{equation}
	\PP{ \rho > y } = \frac{1}{y},  \ \ \ y \geq 1.
	\end{equation}
	The density of  $\rho$ is 
	\begin{equation}
	f _\rho(x) = \frac{1}{x^2},  \ \ \  x\geq 1,
	\end{equation} 
	\begin{equation}
	\E \rho = \int\limits _{1} ^\infty \frac{xdx}{x^2} = \infty,
	\end{equation}
	and for $t \geq 1$
	\begin{equation}
	\E (t \wedge \rho) = t \PP{ \rho > t} + \EE{ \rho \1\{\rho \leq t\} } = t \frac 1t + \int\limits _{1} ^t \frac{xdx}{x^2}
	=  1 + \ln t,
	\end{equation}
	hence
	\begin{equation}
	\int\limits _0 ^\infty
	e^{- \E( \rho \wedge t)} dt
	\geq	\int\limits _1 ^\infty
	e^{- 1 - \ln t} dt
	= e^{-1}
	\int\limits _1 ^\infty
	\frac { dt}{t} = \infty.
	\end{equation}
\end{rmk}

\begin{rmk}
	An alternative proof of (iii) of Theorem \ref{main thm} 
	based on the recurrence properties of Markov processes \cite{MeynTweedie3}
	can be found in the first Arxiv version of this paper.
\end{rmk}



\bibliographystyle{alpha}
\bibliography{Sinus}

\end{document}